\documentclass[a4paper,12pt,reqno]{amsart}
\usepackage{hyperref}

\begin{document}
\title[$q-$difference equations concerning with $q-$gamma function] {$q-$difference equations concerning with $q-$gamma function}

\author[Zhi-Bo Huang,~Ran-Ran Zhang]
{Zhi-Bo Huang$^{*}$,~ Ran-Ran Zhang$^{*}$}

 {\small\thanks{$^{*}$Corresponding author:} }
 \maketitle
 \vskip .2cm\par
\vspace{-0.6cm}
\begin{center}
{\bf ABSTRACT}
\end{center}
We consider a family of solutions of $q-$difference Riccati equation, and prove the meromorphic solutions of $q-$difference Riccati equation and corresponding second order $q-$difference equation are concerning with $q-$gamma function. The growth and value distribution of differences on solutions of $q-$difference Riccati equation are also investigated.
\\\quad\\
{\bf Keywords:}$q-$difference Riccati
equation;$q-$difference equation; $q-$gamma function \\\quad\\
{\bf MSC(2010):} Primary:39B32; Secondary:30D35

\vskip .8cm\par

\maketitle
\numberwithin{equation}{section}
\newtheorem{theorem}{Theorem}[section]
\newtheorem{lemma}[theorem]{Lemma}
\newtheorem{proposition}[theorem]{Proposition}
\newtheorem{remark}[theorem]{Remark}
\newtheorem{definition}[theorem]{Definition}
\allowdisplaybreaks

\section{Introduction}
 \vskip.2cm\par
Let $q\in\mathbb{C}$ be such that $0<|q|<1$.
It is well know that $q$-Gamma function $\Gamma_{q}(x)$ is defined by
\begin{equation*}
  \begin{split}
\Gamma_{q}(x):= \frac{(q; q)_{\infty}}{\left(q^{x};
q\right)_{\infty}}(1-q)^{1-x},
  \end{split}
\end{equation*}
where $(a; q)_{\infty}=\prod_{k=0}^{\infty}\left(1-a q^{k}\right).$
It is a meromorphic function with poles at $x=-n\pm 2\pi i k/\log
q$, where $k$
 and $n$ are non-negative integers, see \cite{AAR}.
  \vskip .2cm\par
By defining
\begin{equation*}
  \begin{split}
 \gamma_{q}(z):=(1-q)^{x-1} \Gamma_{q}(x),\quad z=q^{x},
   \end{split}
\end{equation*}
 and $ \gamma_{q}(0):=(q; q)_{\infty}$, we obtain that $ \gamma_{q}(z)$ is a meromorphic function of zero order with no zeros, having its poles at
 $\left\{q^{k}\right\}_{k=0}^{\infty}$.
  \vskip .2cm\par
  We easily conclude that the first order linear $q$-difference equation
\begin{equation*}
  h(q z)=(1-z)h(z)
\end{equation*}
  is solved by the function $\gamma_{q}(z)$.
 \vskip .2cm\par
 We then consider that a general first order linear $q$-difference
equation
\begin{equation}\label{eq1.1}
  h(q z)=a(z)h(z),
\end{equation}
  where $a(z)$ is a rational function. If $a(z)\equiv a$ is a constant, equation (\ref{eq1.1}) is solvable in terms of rational functions if and only if
  $\log_{q} a$ is an integer.
  If $a(z)$ is nonconstant rational function, let $\alpha_{i}, i=1,2,\cdots, n$ and $\beta_{j}, j=1,2,\cdots, m$ be the zeros and poles of
   $a(z)$, respectively, repeated according to their multiplicities. Then $a(z)$ can be written in the form
    \begin{equation*}
   a(z)=\frac{c(1-z/\alpha_{1})(1-z/\alpha_{2})\cdots(1-z/\alpha_{n})}{(1-z/\beta_{1})(1-z/\beta_{2})\cdots (1-z/\beta_{m})},
 \end{equation*}
   where $c\neq 0$ is a complex number depending on $a(z)$. Thus, equation (\ref{eq1.1}) is solved by
\begin{equation}\label{eq1.2}
   h(z)=z^{\log_{q} c}\frac{\gamma_{q}(z/\alpha_{1})\gamma_{q}(z/\alpha_{2})\cdots \gamma_{q}(z/\alpha_{n})}
   {\gamma_{q}(z/\beta_{1})\gamma_{q}(z/\beta_{2})\cdots\gamma_{q}(z/\beta_{m}) },
   \end{equation}
   which is meromorphic if and only if $\log_{q} c$ is an integer.
 \vskip .2cm\par
In this paper, we are  concerned with the $q$-difference Riccati equation
\begin{equation}\label{eq1.3}
 f(q z)=\frac{A(z)+f(z)}{1-(q-1) z f(z)},
 \end{equation}
 and second order linear $q$-difference equation
 \begin{equation}\label{eq1.4}
 \Delta_{q}^{2} y(z)+\frac{A(z)}{(q-1) z}y(z)=0,
 \end{equation}
where $q\in\mathbb{C}\backslash\{0\}, |q|\neq 1$, $A(z)$ is
meromorphic function. Here, for a meromorphic function $y(z)$, the
$q$-difference operator $\Delta_{q}$ is defined by
$\Delta_{q}y(z)=\frac{y(q z)-y(z)}{(q-1)z}$ and
$\Delta_{q}^{n+1}y(z)=\Delta_{q}\left(\Delta_{q}^{n}y(z)\right),
 n=1,2,\cdots$\cite[pp.488]{AAR}.
\vskip.2cm\par
Throughout this paper, we assume that the reader is familiar with
the fundamental results and the standard notations of Nevanlinna's
value distribution theorem of meromorphic function, see, e.g.,
\cite{Hayman}.
 \vskip .2cm\par
 Recently, a number of papers, see e.g. \cite{Chen1, Chen2, CF, HK1, HK2, HKLRT, Ishizaki, LY1, Wang, ZC}, focused on complex difference equations and difference analogues of Nevanlinna theory.  $q-$difference counterparts are also investigated \cite{BHMK, LY1,
 ZK}. But there are
 only few  papers concerning with the properties of meromorphic solutions of $q-$difference
 equations, see e.g. \cite{Huang, JC, Zheng}.
  \vskip .2cm\par
  The remainder of the paper is organized as follows. A system of solutions of $q-$difference Riccati equation (\ref{eq1.3}) is stated in section 2. Section 3 contains the relationships between $q-$difference equations and $q-$gamma function. The growth and value distribution of differences on solutions of $q-$difference Riccati equation (\ref{eq1.3}) is investigated in Section 4.
    \vskip .5cm\par
  \section{A system of solutions of $q-$difference Riccati equation}
   \vskip.2cm\par
   Let $f_{1}, f_{2}, f_{3}$ be distinct meromorphic solutions of
 differential Riccati equation
 \begin{equation}\label{eq2.1}
w^{'}(z)+w(z)^{2}+A(z)=0,
\end{equation}
Then (\ref{eq2.1}) possesses a one parameter family of meromorphic solutions
$(f_{c})_{c\in\mathbb{C}}$, see e.g., \cite[pp.371-373]{Bank}.
 \vskip.2cm\par
 Ishizaki extended this property to the difference Riccati equation
\begin{equation}\label{eq2.2}
\Delta f(z)+\frac{f(z)^{2}+A(z)}{f(z)-1}=0,
\end{equation}
  and obtained the following
 difference analogue of this property.
 \vskip.2cm\par
 {\bf Theorem 2.A}\cite[Proposition 2.1]{Ishizaki}. Suppose that
(\ref{eq2.2}) possesses three distinct meromorphic solutions
$f_{1}(z), f_{2}(z)$ and $f_{3}(z)$. Then any meromorphic solution
$f(z)$ of (\ref{eq2.2}) can be represented by
\begin{equation}\label{eq2.3}
\begin{split}
f(z)=\frac{f_{1}(z)f_{2}(z)-f_{2}(z)f_{3}(z)-f_{1}(z)f_{2}(z)Q(z)+f_{1}(z)f_{3}(z)Q(z)}
{f_{1}(z)-f_{3}(z)-f_{2}(z)Q(z)+f_{3}(z)Q(z)},
\end{split}
\end{equation}
where $Q(z)$ is a periodic function of period 1. Conversely, if for
any periodic function $Q(z)$ of period 1, we define a function
$f(z)$ by (\ref{eq2.3}), then $f(z)$ is a meromorphic solution of
(\ref{eq2.2}).
 \vskip.2cm\par
  We then extend these properties to the $q-$difference Riccati equation (\ref{eq1.3}) and obtained a
 $q-$difference analogue of these properties as follows.
 \vskip.2cm\par
 {\bf Theorem 2.1}.  Suppose that (\ref{eq1.3}) possesses three distinct meromorphic solutions $f_{1}(z), f_{2}(z)$ and $f_{3}(z)$. Then any
 meromorphic solution $f(z)$ of (\ref{eq1.3}) can be represented by
\begin{equation}\label{eq2.4}
 f(z)=\frac{f_{1}(z)f_{3}(z)-f_{1}(z)f_{2}(z)+f_{1}(z)f_{2}(z)\phi(z)-f_{2}(z)f_{3}(z)\phi(z)}
 {f_{3}(z)-f_{2}(z)-f_{3}(z)\phi(z)+f_{1}(z)\phi(z)},
 \end{equation}
 where $\phi(z)$ is a meromorphic function satisfying $\phi(q z)=\phi(z)$. Conversely, if for any meromorphic function $\phi(z)$
  satisfying $\phi(q z)=\phi(z)$, we define a function $f(z)$ by (\ref{eq2.4}), then $f(z)$ is a meromorphic solution of (\ref{eq1.3}).
 \vskip.2cm\par
 {\bf Proof of Theorem 2.1}. Using a similar proof of Theorem 2 in \cite{Huang}. Let $h_{j}(z), j=1,2,3,4$ be distinct meromorphic functions. We denote a cross ratio of  $h_{j}(z), j=1,2,3,4$ by
\begin{equation*}
R(h_{1}, h_{2}, h_{3}, h_{4}; z):=\frac{h_{4}(z)
-h_{1}(z)}{h_{4}(z)-h_{2}(z)}:\frac{h_{3}(z)-h_{1}(z)}{h_{3}(z)-h_{2}(z)}.
\end{equation*}
  \vskip .2cm\par
 We first show
that $f(z)$, distinct from  $f_{1}(z), f_{2}(z)$ and $f_{3}(z)$, is a meromorphic solution of (\ref{eq1.3}) if and only
if $R(q z)=R(z)$, where $R(z)=R(f_{1}, f_{2}, f_{3}, f; z)$.  Thus, we conclude from (\ref{eq1.3}) that
 \begin{equation*}
 \begin{split}
 R(q z)&=\frac{f(qz)
-f_{1}(qz)}{f(qz)-f_{2}(qz)}:\frac{f_{3}(qz)-f_{1}(qz)}{f_{3}(qz)-f_{2}(qz)}\\
&=\frac{\frac{[(q-1)z A(z)+1][f(z)-f_{1}(z)]}{[1-(q-1)z
f(z)][1-(q-1)z f_{1}(z)]}} {\frac{[(q-1)z
A(z)+1][f(z)-f_{2}(z)]}{[1-(q-1)z f(z)][1-(q-1)z f_{2}(z)]}} :
\frac{\frac{[(q-1)z A(z)+1][f_{3}(z)-f_{1}(z)]}{[1-(q-1)z
f_{3}(z)][1-(q-1)z f_{1}(z)]}}
{\frac{[(q-1)z A(z)+1][f_{3}(z)-f_{2}(z)]}{[1-(q-1)z f_{3}(z)][1-(q-1)z f_{2}(z)]}}\\
 &=\frac{f(z)-f_{1}(z)}{f(z)-f_{2}(z)}:\frac{f_{3}(z)-f_{1}(z)}{f_{3}(z)-f_{2}(z)}=
 R(z).
 \end{split}
  \end{equation*}
     \vskip.2cm\par
   On the other hand, if $R(q z)=R(z)$, then
\begin{equation*}
\begin{split}
   &\frac{f(q z)-\frac{A(z)+f_{1}(z)}{1-(q-1) z f_{1}(z)}}{f(q z)-\frac{A(z)+f_{2}(z)}{1-(q-1) z
   f_{2}(z)}}:
   \frac{\frac{A(z)+f_{3}(z)}{1-(q-1) z f_{3}(z)}-\frac{A(z)+f_{1}(z)}{1-(q-1) z f_{1}(z)}}
   {\frac{A(z)+f_{3}(z)}{1-(q-1) z f_{3}(z)}-\frac{A(z)+f_{2}(z)}{1-(q-1) z f_{2}(z)}}\\
  &=\frac{f(z)-f_{1}(z)}{f(z)-f_{2}(z)}:\frac{f_{3}(z)-f_{1}(z)}{f_{3}(z)-f_{2}(z)},
  \end{split}
 \end{equation*}
 and so,
  \begin{equation}\label{eq2.5}
\begin{split}
   &\frac{f(q z)-\frac{A(z)+f_{1}(z)}{1-(q-1) z f_{1}(z)}}{f(q z)-\frac{A(z)+f_{2}(z)}{1-(q-1) z
   f_{2}(z)}}:
   \frac{\frac{[(q-1)z A(z)-1](f_{3}(z)-f_{1}(z))}{[1-(q-1)z f_{3}(z)][1-(q-1)z f_{1}(z)]}}
   {\frac{[(q-1)z A(z)-1](f_{3}(z)-f_{2}(z))}
  {[1-(q-1)z f_{3}(z)][1-(q-1)z f_{2}(z)]}}\\
  &=\frac{f(z)-f_{1}(z)}{f(z)-f_{2}(z)}:\frac{f_{3}(z)-f_{1}(z)}{f_{3}(z)-f_{2}(z)}.
  \end{split}
 \end{equation}
   \vskip.2cm\par
   We then conclude from (\ref{eq2.5}) that $f(q z)=\frac{A(z)+f(z)}{1-(q-1)
z f(z)}$, which shows that $f(z)$ satisfies (\ref{eq1.3}).
\vskip.2cm\par
 Thus, for any meromorphic function $\phi(z)$
  satisfying $\phi(q z)=\phi(z)$, we define $f(z)$ by
\begin{equation*}
  R(f_{1}, f_{2}, f_{3}, f; z)=\phi(z).
  \end{equation*}
Then $f(z)$ is represented by (\ref{eq2.4}), and also satisfies
(\ref{eq1.3}). The proof of Theorem 2.1 is completed.
 \vskip.2cm\par
 It is difficult
for us to detect the properties of meromorphic solutions since the parameter function $Q(z)$ in Theorem 2.A and
$\phi(z)$ in Theorem 2.1 appear more than one time. Furthermore,
we note that $f(z)\neq f_{2}(z)$ in Theorem 2.1. This shows
that the representation of (\ref{eq2.4}) cannot represent all
meromorphic solutions of $q-$difference Riccati equation
(\ref{eq1.3}). Thus, we can use a new method  used in \cite[Theorem 8.3.4]{Chen3}, and prove a
family of  solutions of $q-$difference Riccati equation
(\ref{eq1.3}).
 \vskip.2cm\par
  {\bf Theorem 2.2.} Let $q\in\mathbb{C}\backslash\{0\}, |q|\neq 1$, and $A(z)$
  be meromorphic function with $A(z)\neq-\frac{1}{(q-1)z}$. If
  $q-$difference Riccati equation (\ref{eq1.3}) possesses three distinct meromorphic solutions $f_{0}(z), f_{1}(z)$ and $f_{2}(z)$,
  then all meromorphic solutions of
   $q-$difference Riccati equation (\ref{eq1.3}) constitute  a one
   parameter family
 \begin{equation}\label{eq2.6}
\left\{f_{0}(z),~f(z)=\frac{(f_{1}(z)-f_{0}(z))(f_{2}(z)-f_{0}(z))}
 {\phi(z)(f_{2}(z)-f_{1}(z))+(f_{2}(z)-f_{0}(z))}+f_{0}(z)\right\},
 \end{equation}
where $\phi(z)$ is any constant in $\mathbb{C}$, or any non-zero
meromorphic function with $\phi(q z)=\phi(z)$, as $\phi(z)\equiv 0,
f(z)=f_{1}(z)$; as $\phi(z)\equiv -1, f(z)=f_{2}(z)$.
 \vskip.2cm\par
 In particular, If $\phi(z)$ is any constant in $\mathbb{C}$, we obtain
  \vskip.2cm\par
 {\bf Corollary 2.1.} Let $q\in\mathbb{C}\backslash\{0\}, |q|\neq 1$, and $A(z)$
  be meromorphic function with $A(z)\neq-\frac{1}{(q-1)z}$. If
  $q-$difference Riccati equation (\ref{eq1.3}) possesses three distinct rational solutions $f_{0}(z), f_{1}(z)$ and $f_{2}(z)$,
  then $q-$difference Riccati equation (\ref{eq1.3}) has infinitely
  many rational solutions.
  \vskip.2cm\par
  We now list some preliminaries to prove Theorem 2.2.
 \vskip.2cm\par
 {\bf Lemma 2.1.} Let $q\in\mathbb{C}\backslash\{0\}, |q|\neq 1$, and $A(z)$
  be meromorphic function with $A(z)\neq-\frac{1}{(q-1)z}$. If $f(z)$
  is a meromorphic solution of $q-$difference Riccati equation
  (\ref{eq1.3}), then
\begin{equation*}
1-(q-1)z f(z)\not\equiv 0~~and ~~ 1+(q-1)z f(qz)\not\equiv 0.
 \end{equation*}
 \vskip.2cm\par
 {\bf Proof of Lemma 2.1.} If $1-(q-1)z f(z)\equiv 0$, then
 $f(z)=\frac{1}{(q-1)z}$. Now substituting $f(z)=\frac{1}{(q-1)z}$
 into (\ref{eq1.3}), and noting that
 $A(z)\neq-\frac{1}{(q-1)z}$, we conclude that
\begin{equation*}
\frac{1}{(q-1)q z}=\frac{A(z)+\frac{1}{(q-1)z}}{1-(q-1)z\cdot
\frac{1}{(q-1)z}}=\frac{(q-1)z A(z)+1}{0}=\infty.
 \end{equation*}
 This yields that $q=0$ or $q=1$, a contradiction.
  \vskip.2cm\par
  If $1+(q-1)z f(qz)\equiv 0$, then $f(q z)=-\frac{1}{(q-1)z}$ and
  $f(z)=-\frac{q}{(q-1)z}$. Now substituting these into
  (\ref{eq1.3}), we deduce that $A(z)=-\frac{1}{(q-1)z}$, a contradiction.
  \vskip.2cm\par
  {\bf Lemma 2.2.} Let $q\in\mathbb{C}\backslash\{0\}, |q|\neq 1$, $A_{1}(z)$ and $A_{0}(z)$ be nonzero meromorphic functions. If $q-$difference equation
 \begin{equation}\label{eq2.7}
A_{1}(z) y(q z)+A_{0}(z)y(z)=0
 \end{equation}
 has a nonzero meromorphic solution $y_{0}(z)$, then all
 meromorphic solutions of (\ref{eq2.7}) constitute a one parameter family
 \begin{equation*}
\{y(z)=\phi(z)y_{0}(z)\},
 \end{equation*}
 where $\phi(z)$ is any constant in $\mathbb{C}$, or any nonzero
meromorphic function with $\phi(q z)=\phi(z)$.
 \vskip.2cm\par
 {\bf Proof of Lemma 2.2.} Since $y_{0}(z)$ is a nonzero
 meromorphic solution of  (\ref{eq2.7}), we easily conclude
 that $y(z)=\phi(z)y_{0}(z)$ is also a meromorphic solution of
 (\ref{eq2.7}) for any constant $\phi(z)$ in $\mathbb{C}$, or any
 non-zero meromorphic function with $\phi(q z)=\phi(z)$.
  \vskip.2cm\par
  On the other hand, if $y(z)$ is also meromorphic solution of
  (\ref{eq2.7}), we conclude from (\ref{eq2.7}) that
  \begin{equation*}
\frac{y(q z)}{y_{0}(q z)}\equiv\frac{y(z)}{y_{0}(z)}.
 \end{equation*}
Set $\phi(z)=\frac{y(z)}{y_{0}(z)}$. Then $\phi(z)$ is a constant in
$\mathbb{C}$, or a nonzero meromorphic function with $\phi(q z)=\phi(z)$.
This shows that $y(z)=\phi(z)y_{0}(z).$
 \vskip.2cm\par
 We now give the proof of Theorem 2.2.
 \vskip.2cm\par
 {\bf Proof of Theorem 2.2.} Since $f_{0}(z), f_{1}(z)$ and $f_{2}(z)$ are three distinct meromorphic solutions
 of $q-$difference Riccati equation (\ref{eq1.3}), we set
 \begin{equation}\label{eq2.8}
u_{j}(z)=\frac{1}{f_{j}(z)-f_{0}(z)}~,~j=1,2.
 \end{equation}
 Obviously, $u_{1}(z)\not\equiv u_{2}(z)$ and
 $f_{j}(z)=\frac{1}{u_{j}(z)}+f_{0}(z),~j=1,2$.
  \vskip.2cm\par
   Now, substituting $f_{j}(z)=\frac{1}{u_{j}(z)}+f_{0}(z),~j=1,2$ into
  (\ref{eq1.3}), and noting that $f_{0}(z)$ is also a meromorphic
  solution of (\ref{eq1.3}), we conclude that
 \begin{equation*}
 [1+(q-1)z f_{0}(q z)]u_{j}(q z)- [1-(q-1)z
 f_{0}(z)]u_{j}(z)+(q-1)z=0.
 \end{equation*}
 Set
 \begin{equation*}
\alpha_{1}(z)= 1+(q-1)z f_{0}(q z)~~and ~~ \alpha_{0}(z)=(q-1)z
f_{0}(z)-1.
 \end{equation*}
 Then we deduce from Lemma 2.1 that $\alpha_{1}(z)\not\equiv 0$ and $\alpha_{0}(z)\not\equiv 0$, and $u_{j}(z),~j=1,2$ are two distinct
  meromorphic solutions of $q-$difference
 equation
\begin{equation}\label{eq2.9}
 \alpha_{1}(z)u(q z)+\alpha_{0}(z)u(z)+(q-1)z=0.
 \end{equation}
 Thus, $u_{0}(z)=u_{1}(z)-u_{2}(z)$
 is a nonzero meromorphic solution of $q-$difference
 equation
\begin{equation}\label{eq2.10}
 \alpha_{1}(z)u(q z)+\alpha_{0}(z)u(z)=0,
 \end{equation}
 which is a corresponding linear homogeneous $q-$difference equation
 of (\ref{eq2.9}).
  \vskip.2cm\par
  Therefore, we deduce from Lemma 2.2 that all
 meromorphic solutions of (\ref{eq2.10}) constitute a one parameter family
 \begin{equation*}
H(y(z))=\{y(z)=\phi(z)u_{0}(z)\},
 \end{equation*}
 where $\phi(z)$ is any constant in $\mathbb{C}$, or any non-zero
meromorphic function with $\phi(q z)=\phi(z)$. This yields that
$q-$difference equation (\ref{eq2.9}) has a general solution
 \begin{equation}\label{eq2.11}
 \begin{split}
 u(z)&=y(z)+u_{1}(z)=Q(z)u_{0}(z)+u_{1}(z)\\
 &=Q(z)[u_{1}(z)-u_{2}(z)]+u_{1}(z)\\
 &=\frac{Q(z)[f_{2}(z)-f_{1}(z)]}{[f_{1}(z)-f_{0}(z)][f_{2}(z)-f_{0}(z)]}+\frac{1}{f_{1}(z)-f_{0}(z)}.
 \end{split}
 \end{equation}
\vskip.2cm\par
 We now suppose that $f(z)(\not\equiv f_{0}(z))$ is a meromorphic
 solution of (\ref{eq1.3}), and conclude from the argumentation of
 (\ref{eq2.9}) that $u(z)=\frac{1}{f(z)-f_{0}(z)}$ is also a meromorphic
 solution of (\ref{eq2.9}). Thus, we deduce from (\ref{eq2.11}) that
 there exists a constant $\phi(z)$ in $\mathbb{C}$, or any non-zero
meromorphic function $\phi(z)$ with $\phi(q z)=\phi(z)$ such that
  \begin{equation*}
\frac{1}{f(z)-f_{0}(z)}=\frac{\phi(z)[f_{2}(z)-f_{1}(z)]}{[f_{1}(z)-f_{0}(z)][f_{2}(z)-f_{0}(z)]}+\frac{1}{f_{1}(z)-f_{0}(z)}.
 \end{equation*}
 Therefore, we obtain that
  \begin{equation}\label{eq2.12}
 f(z)=\frac{[f_{1}(z)-f_{0}(z)][f_{2}(z)-f_{0}(z)]}{\phi(z)[f_{2}(z)-f_{1}(z)]+[f_{2}(z)-f_{0}(z)]}+f_{0}(z),
 \end{equation}
 where $\phi(z)$ is any constant in $\mathbb{C}$, or any non-zero
meromorphic function with $\phi(q z)=\phi(z)$.
 This shows that any meromorphic solution $f(z)(\not\equiv
 f_{0}(z))$ of (\ref{eq1.3}) has the form (\ref{eq2.12}).
 \vskip.2cm\par
  We then affirm that any meromorphic function $f(z)(\not\equiv f_{0}(z))$ denoted by
 (\ref{eq2.12}) must be a meromorphic solution of (\ref{eq1.3}). In
 fact, we conclude from (\ref{eq2.11}) and (\ref{eq2.12}) that
 \begin{equation}\label{eq2.13}
 f(z)=\frac{1}{u(z)}+f_{0}(z),
 \end{equation}
 where $u(z)$ satisfies the $q-$difference equation (\ref{eq2.9}).
 Thus, we further conclude from (\ref{eq2.9}),(\ref{eq1.3}) and the assumption that $f_{0}(
 z)$ is a meromorphic solution of (\ref{eq1.3}), that
 \begin{equation}\label{eq2.14}
  \begin{split}
 f(q z)&=\frac{1}{u(q z)}+f_{0}(q z)\\
 &=\frac{\alpha_{1}(z)}{-\alpha_{0}(z)u(z)-(q-1)z}+f_{0}(q
 z)\\
 &=\frac{1+(q-1)z f_{0}(q z)}{[1-(q-1)z f_{0}(z)]u(z)-(q-1)z}+f_{0}(q
 z)\\
 &=\frac{1+[1-(q-1)z f_{0}(z)]u(z)f_{0}(q z)}{[1-(q-1)z f_{0}(z)]u(z)-(q-1)z}\\
 &=\frac{1+[1-(q-1)z f_{0}(z)]u(z)\cdot \frac{A(z)+f_{0}(
 z)}{1-(q-1)z f_{0}(
 z)}}{[1-(q-1)z f_{0}(z)]u(z)-(q-1)z}\\
 &=\frac{1+[A(z)+f_{0}(
 z)]u(z)}{[1-(q-1)z f_{0}(z)]u(z)-(q-1)z}.
  \end{split}
 \end{equation}
 \vskip.2cm\par
 On the other hand, we can obtain from  (\ref{eq2.13}) that
 \begin{equation}\label{eq2.15}
   \begin{split}
 \frac{A(z)+f(z)}{1-(q-1)z f(z)}&=\frac{A(z)+\frac{1}{u(z)}+f_{0}(z)}{1-(q-1)z
 \left[\frac{1}{u(z)}+f_{0}(z)\right]}\\
 &=\frac{1+[A(z)+f_{0}(
 z)]u(z)}{[1-(q-1)z f_{0}(z)]u(z)-(q-1)z}.
   \end{split}
 \end{equation}
 \vskip.2cm\par
  Therefore, we deduce from (\ref{eq2.14}) and (\ref{eq2.15}) that
 \begin{equation*}
 f(q z)=\frac{A(z)+f(z)}{1-(q-1)z f(z)},
 \end{equation*}
 which shows that meromorphic function $f(z)(\not\equiv f_{0}(z))$ denoted by
 (\ref{eq2.12}) is a meromorphic solution of (\ref{eq1.3}). The
 proof of Theorem 2.2 is completed.
  \vskip.5cm\par
  If $q-$difference Riccati equation
  (\ref{eq1.3}) possesses a rational solution $f_{0}(z)$ such that
  $f_{0}(z)\not\rightarrow 0$ as $z\rightarrow \infty$, we further obtain
  \vskip.2cm\par
   {\bf Theorem 2.3.} Let $q\in\mathbb{C}\backslash\{0\}, |q|\neq 1$, and $A(z)$
  be meromorphic function. If
  $q-$difference Riccati equation (\ref{eq1.3}) possesses a rational  solution $f_{0}(z)$ such that
  $f_{0}(z)\not\rightarrow 0$ as $z\rightarrow \infty$,
  then
   $q-$difference Riccati equation (\ref{eq1.3}) has at most two
   rational solutions.
  \vskip.2cm\par
  {\bf Proof of Theorem 2.3.} Contrary to the assumption, we suppose
  that (\ref{eq1.3}) has three distinct rational solutions $f_{0}(z), f_{1}(z)$ and $f_{2}(z)$, where  $f_{0}(z)\not\rightarrow 0$ as $z\rightarrow \infty$.
    \vskip.2cm\par
  Set
 \begin{equation*}
u_{j}(z)=\frac{1}{f_{j}(z)-f_{0}(z)}~,~j=1,2.
 \end{equation*}
 Obviously, $u_{1}(z)$ and $u_{2}(z)$ are rational functions with $u_{1}(z)\not\equiv u_{2}(z)$, and
 $f_{j}(z)=\frac{1}{u_{j}(z)}+f_{0}(z),~j=1,2$.
  \vskip.2cm\par
  Now, substituting $f_{j}(z)=\frac{1}{u_{j}(z)}+f_{0}(z),~j=1,2$ into
  (\ref{eq1.3}), and noting that $f_{0}(z)$ is also a rational
  solution of (\ref{eq1.3}), we conclude that
 \begin{equation*}
 [1+(q-1)z f_{0}(q z)]u_{j}(q z)- [1-(q-1)z
 f_{0}(z)]u_{j}(z)+(q-1)z=0 ,~j=1,2.
 \end{equation*}
 This shows that  $u_{1}(z)$ and $u_{2}(z)$ are two distinct rational
 solutions
 of $q-$difference equation
\begin{equation}\label{eq2.16}
 [1+(q-1)z f_{0}(q z)]u(q z)- [1-(q-1)z
 f_{0}(z)]u(z)+(q-1)z=0,
 \end{equation}
 and $u_{0}(z)=u_{1}(z)-u_{2}(z)$
 is a nonzero rational solution of
\begin{equation}\label{eq2.17}
 [1+(q-1)z f_{0}(q z)]u(q z)- [1-(q-1)z
 f_{0}(z)]u(z)=0.
 \end{equation}
\vskip.2cm\par
 Since $f_{0}(z)(\not\rightarrow 0, z\rightarrow\infty)$ and $u_{0}(z)$ are both rational functions, we can
 set
\begin{equation*}
 f_{0}(z)=\frac{P(z)}{Q(z)} ~and ~u_{0}(z)=\frac{U(z)}{V(z)},
 \end{equation*}
 where $P(z), Q(z), U(z)$ and $V(z)$ are nonzero polynomials with
 $\deg P(z)\geq \deg Q(z)$.
 \vskip.2cm\par
 Now substituting $f_{0}(z)=\frac{P(z)}{Q(z)}$ and $u_{0}(z)=\frac{U(z)}{V(z)}$ into
 (\ref{eq2.17}), we conclude that
 \begin{equation}\label{eq2.18}
    \begin{split}
 &Q(z)Q(q z)U(q z)V(z)+(q-1)z P(q z)Q(z)U(q z)V(z)\\
&-Q(z)Q(q z)U(z)V(q z)+(q-1)z P(z)Q(q z)U(z)V(q z)=0.
   \end{split}
 \end{equation}
 \vskip.2cm\par
 We can obtain that
 \begin{equation*}
    \begin{split}
&\deg\{Q(z)Q(q z)U(q z)V(z)\}
=\deg\{Q(z)Q(q z)U(z)V(q z)\}\\
&<\deg\{(q-1)z P(q z)Q(z)U(q z)V(z)\}
=\deg\{(q-1)z P(z)Q(q z)U(z)V(q z)\},
 \end{split}
 \end{equation*}
 and at most one of  the coefficient of power $z^{\deg\{Q(z)Q(q z)U(q
 z)V(z)\}}$ and the coefficient of power $z^{\deg\{(q-1)z P(q z)Q(z)U(q
 z)V(z)\}}$ is zero. These all show that the degree of left hand
 side of (\ref{eq2.18}) is great than 1, and yield a contradiction.
 \vskip.2cm\par
 Thus, (\ref{eq1.3}) has at most two rational solutions. The proof of Theorem 2.3 is completed.
  \vskip.2cm\par
  We now present two examples to show that Theorem 2.3 remain valid.
 \vskip.2cm\par
  {\bf Example 2.1.} Let $q=\frac{1}{2}$. Then ration function $f_{0}(z)=2z+4$ solves the $q-$difference
  Riccati equation
 \begin{equation}\label{eq2.19}
 f\left(\frac{1}{2}z\right)=\frac{z^{3}+6 z^{2}+7
  z+f(z)}{1+\frac{z}{2}f(z)}
 \end{equation}
 of type (\ref{eq1.3}), and $f_{0}(z)=2z+4 \rightarrow \infty$ as $z\rightarrow\infty$.  Suppose that $f_{1}(z)(\not\equiv f_{0}(z))$ is
 another rational solution of (\ref{eq2.19}). Set
 $u_{1}(z)=\frac{1}{f_{1}(z)-f_{0}(z)}$. Then we conclude that $u_{1}(z)$
 satisfies the $q-$difference equation
\begin{equation}\label{eq2.20}
 (z^{2}+4 z-2)u\left(\frac{z}{2}\right)+2(z^{2}+2z+1)u(z)+z=0.
 \end{equation}
 \vskip.2cm\par
  According to the proof of Theorem 2.2, we note that all
 meromorphic solutions $f(z)$(except exceptional solution
 $f_{0}(z)$) of (\ref{eq2.19}) and all solutions
 $u(z)=\frac{1}{f(z)-f_{0}(z)}$ of (\ref{eq2.20}) are one-one
 corresponding.
 \vskip.2cm\par
 However, the $q-$difference equation, which is the corresponding
 homogeneous difference equation of (\ref{eq2.20}),
 \begin{equation}\label{eq2.21}
 (z^{2}+4 z-2)u\left(\frac{z}{2}\right)+2(z^{2}+2z+1)u(z)=0,
 \end{equation}
 has no nonzero rational solution. Otherwise, suppose that
 $u(z)=\frac{P(z)}{Q(z)}$ is a nonzero rational solution of
 (\ref{eq2.21}), where $P(z)$ and $Q(z)$ are nonzero polynomials
 with degree $\deg P(z)=p$ and $\deg Q(z)=q$ respectively. Then we
 conclude from (\ref{eq2.21}) that
  \begin{equation}\label{eq2.22}
 (z^{2}+4
 z-2)P\left(\frac{z}{2}\right)Q(z)+2(z^{2}+2z+1)P(z)Q\left(\frac{z}{2}\right)=0.
 \end{equation}
 We can easily deduce  that the degree of left
 hand side of (\ref{eq2.22}) is great than 2 since $P(z)$ and  $Q(z)$ are nonzero polynomials, and yields a contradiction. Hence, we obtain that
 (\ref{eq2.19}) has at most two rational solutions $f_{0}(z)=2z+4$
 and $f_{1}(z)=\frac{1}{u_{1}(z)}+2z+4$.
 \vskip.2cm\par
  {\bf Example 2.2.} Let $q=\frac{1}{2}$ and $A(z)=\frac{2(z+1)(z+2)}{z(z^{2}-3z-2)}$. Then ration function $f_{0}(z)=\frac{z-1}{z+1}$ solves the $q-$difference
  Riccati equation
 \begin{equation}\label{eq2.23}
 f\left(\frac{1}{2}z\right)=\frac{A(z)+f(z)}{1+\frac{z}{2}f(z)}
 \end{equation}
 of type (\ref{eq1.3}), and $f_{0}(z)=\frac{z-1}{z+1} \rightarrow 1$ as $z\rightarrow\infty$.  Suppose that $f_{1}(z)(\not\equiv f_{0}(z))$ is
 another rational solution of (\ref{eq2.23}). Set
 $u_{1}(z)=\frac{1}{f_{1}(z)-f_{0}(z)}$. Then we conclude that $u_{1}(z)$
 satisfies the $q-$difference equation
\begin{equation}\label{eq2.24}
 (z^{3}-3 z^{2}-8 z-4)u\left(\frac{z}{2}\right)+(z^{3} +3 z^{2}+4z+4)u(z)+z(z^{2}+3z +2)=0.
 \end{equation}
 \vskip.2cm\par
 By using similar calculation of Example 2.1, the $q-$difference equation, which is the corresponding
 homogeneous difference equation of (\ref{eq2.24}),
 \begin{equation}\label{eq2.25}
(z^{3}-3 z^{2}-8 z-4)u\left(\frac{z}{2}\right)+(z^{3} +3
z^{2}+4z+4)u(z)=0,
 \end{equation}
 has no nonzero rational solution. Thus, we obtain that
 (\ref{eq2.23}) has at most two rational solutions $f_{0}(z)=\frac{z-1}{z+1}$
 and $f_{1}(z)=\frac{1}{u_{1}(z)}+\frac{z-1}{z+1}$.

     \vskip .5cm\par
  \section{  Relationships between $q-$difference  equation and $q-$gamma function}
   \vskip.2cm\par
   In this section, we focus on the relationships between $q-$difference  equation and $q-$gamma function, and firstly obtain the following result.
  \vskip.2cm\par
  {\bf Theorem 3.1}. Let $q\in\mathbb{C}$ with
$0<|q|<1$.
  Suppose that $q$-difference Riccati equation (\ref{eq1.3}) possesses two distinct rational solutions $f_{1}(z)$ and $f_{2}(z)$.
  Then  all meromorphic solutions of $q-$difference Riccati equation (\ref{eq1.3}) are concerned with $q$-gamma
  function.
\vskip .2cm\par
{\bf Proof of Theorem 3.1}.  Since $f_{1}(z)$ and $f_{2}(z)$ are two
distinct rational solutions of (\ref{eq1.3}), we construct a
M\"{o}bius translation
\begin{equation}\label{eq3.1}
 \begin{split}
f(z)=\frac{f_{1}(z) h(z)+f_{2}(z)}{h(z)+1}.
 \end{split}
\end{equation}
Then $\sigma(h)=\sigma(f)=0$. Substituting (\ref{eq3.1}) into
(\ref{eq1.3}), we conclude that
\begin{equation}\label{eq3.2}
 \begin{split}
h(q z)=\frac{1-(q-1)z f_{1}(z)}{1-(q-1)z f_{2}(z)}h(z),
 \end{split}
\end{equation}
which is type of (\ref{eq1.1}). Thus, meromorphic solution $h(z)$ of
(\ref{eq3.2}) has the form (\ref{eq1.2}), which is concerned with
$q-$gamma function. The proof of Theorem 3.1 is completed.
\vskip .2cm\par
 Now, we give an example to give a presentation for Theorem 3.1.
 \vskip .2cm\par
{\bf Example 3.1}. Let $q=-\frac{1}{2}$, $A(z)=-\frac{6
z}{(z+1)(z-2)}$ in (\ref{eq1.3}). Then functions
\begin{equation}\label{eq3.3}
 \begin{split}
f_{1}(z)=\frac{1}{z+1}~~ and ~~ f_{2}(z)=\frac{-2}{z+1}
 \end{split}
\end{equation}
satisfy the $q$-difference Riccati equation  (\ref{eq1.3}). Then, by
using the  transformation (\ref{eq3.1}), we can switch
$q-$difference Riccati equation (\ref{eq1.3}) into the type
(\ref{eq3.2}) and conclude that
\begin{equation*}
 \begin{split}
h\left(-\frac{1}{2}
z\right)=\frac{\left(1-\frac{z}{-\frac{2}{3}}\right)}{\left(1-\frac{z}{\frac{1}{2}}\right)}h(z),
 \end{split}
\end{equation*}
which is type of (\ref{eq1.1}), and so
\begin{equation}\label{eq3.4}
 \begin{split}
h(z)=\frac{\gamma_{-\frac{1}{2}}\left(\frac{z}{-\frac{2}{3}}\right)}{\gamma_{-\frac{1}{2}}\left(\frac{z}{\frac{1}{2}}\right)}
=\frac{\gamma_{-\frac{1}{2}}\left(-\frac{3z}{2}\right)}{\gamma_{-\frac{1}{2}}\left(2z\right)}.
 \end{split}
\end{equation}
 \vskip .2cm\par
We then conclude from (\ref{eq3.1}) and (\ref{eq3.4}) that
\begin{equation*}
 \begin{split}
f(z)=\frac{\gamma_{-\frac{1}{2}}\left(-\frac{3z}{2}\right)-2\gamma_{-\frac{1}{2}}\left(2z\right)}
{(z+1)\left(\gamma_{-\frac{1}{2}}\left(-\frac{3z}{2}\right)+\gamma_{-\frac{1}{2}}\left(2z\right)\right)},
 \end{split}
\end{equation*}
which is concerned with $q-$gamma function.
 \vskip .2cm\par
 We second show that solutions of second order $q-$difference equation are also concerning with $q-$gamma function. Thus, we investigate the passage between $q$-difference Riccati
 equation (\ref{eq1.3}) and  second order  $q$-difference
 equation (\ref{eq1.4}), and obtain the following result.
 \vskip .2cm\par
 {\bf Theorem 3.2.} The passage between $q$-difference Riccati
 equation (\ref{eq1.3}) and  second order  $q$-difference
 equation (\ref{eq1.4}) is
 \begin{equation}\label{eq3.5}
 f(z)=-\frac{\Delta_{q}y(z)}{y(z)}=-\frac{y(q z)-y(z)}{(q-1)z y(z)}.
\end{equation}
 \vskip .2cm\par
  {\bf Proof of Theorem 3.2.}  We first prove that $f(z)$ defined as (\ref{eq3.5}) is a
 meromorphic solution of (\ref{eq1.3}) if $y(z)$ is a nontrivial meromorphic solution of (\ref{eq1.4}). In fact, we conclude from (\ref{eq3.5})
  that
\begin{equation}\label{eq3.6}
\begin{split}
 \Delta_{q}^{2} y(z)&=\Delta_{q}(\Delta_{q}y(z))=\Delta_{q}(-f(z)y(z))\\
 &=\frac{-f(q z)y(q z)+f(z)y(z)}{(q-1)z}\\
 &= \frac{-f(q z)[y(z)-(q-1)z f(z) y(z)]+f(z)y(z)}{(q-1)z}\\
 &=\frac{-f(q z)y(z)[1-(q-1)z f(z)]+f(z)y(z)}{(q-1)z}.
 \end{split}
\end{equation}
\vskip.2cm\par
 Thus, we deduce from (\ref{eq1.4}) and (\ref{eq3.6}) that
\begin{equation*}
-f(q z)y(z)[1-(q-1)z f(z)]+f(z)y(z)=-A(z)y(z),
 \end{equation*}
 which implies the desired form of equation (\ref{eq1.3}).
  \vskip .2cm\par
  We second prove that a meromorphic function $y(z)$ satisfying
  (\ref{eq3.5}) is a meromorphic solution of (\ref{eq1.4}) if $f(z)$
  defined as (\ref{eq3.5})
  is a meromorphic solution of (\ref{eq1.3}).
  \vskip .2cm\par
 In fact, we conclude from (\ref{eq3.6}) and (\ref{eq1.3}) that
\begin{equation*}
 \begin{split}
\Delta_{q}^{2} y(z)&=\frac{-f(q z)y(z)[1-(q-1)z f(z)]+f(z)y(z)}{(q-1)z}\\
 &=\frac{-\frac{A(z)+f(z)}{1-(q-1) z f(z)}y(z)[1-(q-1)z f(z)]+f(z)y(z)}{(q-1)z}\\
 &=-\frac{A(z)}{(q-1)z}y(z),
 \end{split}
\end{equation*}
which implies  the desired form of (\ref{eq1.4}).
  \vskip .2cm\par
  Thus, we deduce from Theorem 3.1 and 3.2 that
  \vskip.2cm\par
  {\bf Theorem 3.3}. Let $q\in\mathbb{C}$ with
$0<|q|<1$.
  Suppose that  $q$-difference Riccati equation (\ref{eq1.3}) possesses two distinct rational solutions $f_{1}(z)$ and $f_{2}(z)$.
  Then  all meromorphic solutions of second order $q-$difference equation (\ref{eq1.4})  are concerned with $q$-gamma
  function.
  \vskip .5cm\par
  \section{Value distribution of solutions of $q-$difference Riccati equations}
\vskip .2cm\par
 If $g(z)$ is a transcendental meromorphic solution
of equation
\begin{equation}\label{eq4.1}
g(q z)=R(z, g(z)),
\end{equation}
 where $q\in\mathbb{C}, |q|>1$, and the coefficients of $R(z, g(z))$ are small functions relative to $g(z)$,
  Gundersen et al. \cite{GHLRY} showed that the order
of growth of equation  (\ref{eq4.1})
 is equal to $\log \deg_{g}(R)/ \log |q|$,  where
$\deg_{g}(R)$ is the degree of irreducible rational function $R(z,
g(z))$ in $g(z)$, which means that all transcendental meromorphic
solutions of $q-$difference Riccati equation
  (\ref{eq1.3}) have zero order when $q\in\mathbb{C}, |q|>1$.
  \vskip.2cm\par
  On the
  other hand, second order $q-$difference equation (\ref{eq1.4}) is equivalent to the second order linear
$q-$difference equation
  \begin{equation}\label{eq4.2}
 y\left(q^{2} z\right)-(q+1)y(q z)+q[1+(q-1)z A(z)]y(z)=0.
  \end{equation}
  \vskip.2cm\par
  Bergweiler et al. \cite{BIY} pointed out that all
transcendental meromorphic solutions of equation (\ref{eq4.2})
satisfy $T(r, f)=O((\log r)^2)$ if $q\in\mathbb{C}$ and $0<|q|<1$.
This indicates that all transcendental meromorphic solutions of
equation (\ref{eq1.4}) satisfy $T(r, f)=O((\log r)^2)$ if
$q\in\mathbb{C}$ and $0<|q|<1$. Since (\ref{eq3.5})
 is a passage between (\ref{eq1.3}) and (\ref{eq1.4}),
we deduce that all transcendental meromorphic solutions of
$q-$difference Riccati equation (\ref{eq1.3}) are of zero order if
$q\in\mathbb{C}$ and $0<|q|<1$. Thus, we obtain the following
result.
 \vskip.2cm\par
   {\bf Theorem 4.1.} All transcendental meromorphic solutions of $q-$difference Riccati equation (\ref{eq1.3})
 are of zero order for all $q\in\mathbb{C}\backslash\{0\}$ and $|q|\neq 1$.
 \vskip .2cm\par
Berweiler et. al.\cite{BL} first investigated the existence of
 zeros of $\Delta_{c}f(z)=f(z+c)-f(z)$ and
 $\frac{\Delta_{c}f(z)}{f(z)}$, where
 $c\in\mathbb{C}\backslash\{0\}$, and obtained many profound
 results(see \cite[Theorem 1.2$-$Theorem 1.4 and Theorem 1.6]{BL}).
 Chen and Shon\cite{Chen} then  extended these results of
 \cite{BL} and proved a number of results concerned with the
 existence of zeros and fixed points of $\Delta_{c}f(z)=f(z+c)-f(z)$ and
 $\frac{\Delta_{c}f(z)}{f(z)}$, where
 $c\in\mathbb{C}\backslash\{0\}$(see \cite[Theorem 1$-$Theorem
 6]{Chen}). Zhang and Chen\cite{ZC} further considered the difference Riccati
 equation
 \begin{equation}\label{eq4.3}
f(z+1)=\frac{p(z+1)w(z)+q(z)}{w(z)+p(z)},
\end{equation}
where $p(z)$ and $q(z)$ are small functions relative to $f(z)$, and
obtained the following results.
 \vskip.2cm\par
 {\bf Theorem 4.A}\cite[Theorem 1.1]{ZC}. Let $p(z), q(z)$ be
 meromorphic functions of finite order, and let
 $[p(z+1)f(z)+q(z)]/[f(z)+p(z)]$ be an irreducible function in
 $f(z)$. Suppose that $f(z)$ is an admissible finite order
 meromorphic solution of (\ref{eq4.3}). Set $\Delta f(z)=f(z+1)-f(z)$. Then
  \vskip.2cm\par
   (i)~~$\lambda\left(\frac{1}{\Delta f(z)}\right)=\sigma(\Delta
   f(z))=\sigma(f),~~\lambda\left(\frac{1}{\Delta
   f(z)/f(z)}\right)=\sigma\left(\frac{\Delta
   f(z)}{f(z)}\right)=\sigma(f)$;
   \vskip.2cm\par
   (ii)~~If $q(z)\not\equiv 0$, then
   $\lambda\left(\frac{1}{f}\right)=\lambda(f)=\sigma(f)$;
   \vskip.2cm\par
   (iii)~~If $p(z)\equiv p$ is a constant and $q(z)\equiv s(z)^{2}$,
   where $s(z)$ is a non-constant rational function, then $f(z)$ has
   no Borel exceptional value and $\lambda(\Delta f(z))=\lambda\left(\frac{\Delta
   f(z)}{f}\right)=\sigma(f)$.
   \vskip.2cm\par
    We now consider the value distribution of differences of
   transcendental meromorphic solutions of $q-$difference Riccati
   equation (\ref{eq1.3}) as follows.
   \vskip.2cm\par
   {\bf Theorem 4.2.} Let $A(z)$ be a non-constant rational function,  $q\in\mathbb{C}\backslash\{0\}$ and $|q|\neq 1$.
  Suppose that $f(z)$ is a transcendental meromorphic solution of
  $q-$difference Riccati equation (\ref{eq1.3}). Set $\Delta_{q}f(z)=\frac{f(q
  z)-f(z)}{(q-1)z}$.
  \vskip.2cm\par
  {\bf (1)}  If $A(z)\neq -\frac{1}{(q-1)z}$, then $\Delta_{q}f(z)$ and
  $\frac{\Delta_{q}f(z)}{f(z)}$ have infinitely many poles;
 \vskip.2cm\par
 {\bf (2)}
  If $A(z)=-(q-1)z
  s(z)^{2}$, where $s(z)$ is a non-constant rational function, then $\Delta_{q}f(z)=\frac{f(q z)-f(z)}{(q-1)z}$ and $\frac{\Delta_{q}f(z)}{f(z)}$
  have infinitely many zeros.
    \vskip.2cm\par
     {\bf Remark 4.1} The similar result of Theorem 4.2(2) has been
    obtain in \cite{JC}. For the completeness, we list it again.
         \vskip.2cm\par
     In order to prove to Theorem 4.2, we need some Lemmas.
      \vskip.2cm\par
 {\bf Lemma 4.1}\cite[Theorem 1.2]{BHMK}. Let $f(z)$ be a
 non-constant zero-order meromorphic function, and $q\in
 \mathbb{C}\backslash\{0\}$. Then
 \begin{equation*}
m\left(r, \frac{f(q z)}{f(z)}\right)=o(T(r, f))
\end{equation*}
on a set of logarithmic density 1.
       \vskip.2cm\par
{\bf Lemma 4.2}\cite[Theorem 2.5]{LY1}. Let $f(z)$ be a
 transcendental meromorphic solution of order zero of a
 $q-$difference equation of the form
\begin{equation*}
U_{q}(z,f)P_{q}(z,f)=Q_{q}(z,f),
\end{equation*}
where $U_{q}(z,f)P_{q}(z,f)$ and $Q_{q}(z,f)$ are $q-$difference
polynomials such that the total degree $\deg U_{q}(z,g)=n$ in $f(z)$
and its $q-$shifts, $\deg Q_{q}(z,g)\leq n$. Moreover, we assume
that $U_{q}(z,g)$ contains  just one term of maximal total degree in
$f(z)$ and its $q-$shifts. Then
\begin{equation*}
m(r, P_{q}(z,f))=o(T(r, f))
\end{equation*}
on a set of logarithmic density 1.
\vskip .2cm\par
 {\bf Lemma 4.3}\cite[Theorem 2.2]{BHMK}. Let $f(z)$ be a
nonconstant zero order meromorphic solution of
\begin{equation*}
P(z, f)=0,
\end{equation*}
where $P(z, f)$ is a $q$-difference polynomials in $f(z)$. If $P(z,
\alpha)\not\equiv 0$ for a small function $\alpha(z)$ relative to
$f(z)$, then
\begin{equation*}
m\left(r, \frac{1}{f-\alpha}\right)=o(T(r, f))
\end{equation*}
on a set of logarithmic density 1.
\vskip.2cm\par
{\bf Proof of Theorem 4.2.} {\bf (1)} Since $f(z)$ is a
transcendental meromorphic solution of
  $q-$difference Riccati equation (\ref{eq1.3}), we deduce from
  (\ref{eq1.3}) that
  \begin{equation}\label{eq4.4}
(q-1)z f(q z)f(z)=f(q z)-f(z)-A(z).
\end{equation}
Thus, we conclude from Lemma 4.1, Lemma 4.2 and (\ref{eq4.4}) that
 \begin{equation*}
m(r, f)\leq m\left(r, \frac{f(z)}{f(q z)}\right)-m(r, f(q z))+o(T(r,
f))=o(T(r, f))
\end{equation*}
on a set of logarithmic density 1, and so
\begin{equation}\label{eq4.5}
N(r,f)= T(r, f)+o(T(r, f)),
\end{equation}
 \begin{equation}\label{eq4.6}
m(r, \Delta_{q} f(z))\leq m\left(r, \frac{\Delta_{q}
f(z)}{f(z)}\right)+m(r, f(z))+o(T(r, f))=o(T(r, f))
\end{equation}
on a set of logarithmic density 1.
    \vskip.2cm\par
 We further conclude from (\ref{eq1.3}) that
\begin{equation}\label{eq4.7}
\Delta_{q}f(z)=\frac{1}{(q-1)z}\cdot \frac{A(z)+(q-1)z
f(z)^{2}}{1-(q-1)z f(z)}.
\end{equation}
Thus, we apply Valiron-Mohon'ko Theorem to (\ref{eq4.7}) that
\begin{equation}\label{eq4.8}
T(r, \Delta_{q}f(z))=2 T(r, f)+o(T(r, f)).
\end{equation}
We then obtain from (\ref{eq4.6}) and (\ref{eq4.8}) that
\begin{equation*}
N(r, \Delta_{q}f(z))=2 T(r, f)+o(T(r, f))
\end{equation*}
on a set of logarithmic density 1, and so $\Delta_{q}f(z)$ has
infinitely many poles.
    \vskip.2cm\par
  We note that
    \begin{equation*}
N\left(r, \frac{\Delta_{q}f(z)}{f(z)}\right)\geq N(r,
\Delta_{q}f(z))-N(r, f)= T(r, f)+o(T(r, f))
\end{equation*}
on a set of logarithmic density 1, and so
$\frac{\Delta_{q}f(z)}{f(z)}$ has infinitely many poles.
\vskip.2cm\par
 {\bf (2)} Let
 \begin{equation}\label{eq4.9}
P(z, f)=(q-1)z f(z)f(qz)-f(qz)+f(z)+A(z).
\end{equation}
We then affirm that  $P(z, s(z))\not\equiv 0$ or $P(z,
-s(z))\not\equiv 0$. Otherwise, if $P(z, s(z))\equiv$ and  $P(z,
-s(z))\equiv 0$, we can obtain from (\ref{eq4.9}) that
\begin{equation*}
s(q z)=s(z).
\end{equation*}
This is impossible since $s(z)$ is a non-constant rational function.
Without loss of generality, we assume that  $P(z, s(z))\not\equiv
0$. Thus, we obtain from Lemma 4.3, (\ref{eq1.3}) and (\ref{eq4.9})
that
\begin{equation*}
m\left(r, \frac{1}{f(z)-s(z)}\right)=o(T(r, f))
\end{equation*}
on a set of logarithmic density 1, and so
 \begin{equation}\label{eq4.10}
N\left(r, \frac{1}{f(z)-s(z)}\right)=T(r, f)+o(T(r, f))
\end{equation}
on a set of logarithmic density 1.
    \vskip.2cm\par
 Since $A(z)=-(q-1)z s(z)^{2}$, we can conclude from
    (\ref{eq1.3}) and (\ref{eq4.7}) that
\begin{equation}\label{eq4.11}
\Delta_{q}f(z)=\frac{[f(z)+s(z)][f(z)-s(z)]}{1-(q-1)z f(z)}.
\end{equation}
If $f(z_{0})-s(z_{0})=1-(q-1)z_{0} f(z_{0})=0$, then
$(q-1)z_{0}s(z_{0})=0$. If $f(z_{0})-s(z_{0})=0$ and
$f(z_{0})+s(z_{0})=\infty$, then $s(z_{0})=\infty$. Thus, we deduce
from (\ref{eq4.10}) and (\ref{eq4.11}) that
 \begin{eqnarray*}
N\left(r, \frac{1}{\Delta_{q}f(z)}\right)&=&N\left(r, \frac{1-(q-1)z
f(z)}{[f(z)+s(z)][f(z)-s(z)]}\right) \\
&=&N\left(r, \frac{1}{f(z)-s(z)}\right)=T(r, f)+o(T(r, f))
\end{eqnarray*}
on a set of logarithmic density 1. This shows that
$\frac{1}{\Delta_{q}f(z)}$ has infinitely many zeros.
 \vskip.2cm\par
  We now obtain from (\ref{eq4.11}) that
\begin{equation}\label{eq4.12}
\frac{\Delta_{q}f(z)}{f(z)}=\frac{[f(z)+s(z)][f(z)-s(z)]}{[1-(q-1)z
f(z)]f(z)}.
\end{equation}
By combining (\ref{eq4.5}) and (\ref{eq4.12}), and using similar
method that $\frac{1}{\Delta_{q}f(z)}$ has infinitely many zeros, we
can conclude that $\frac{\Delta_{q}f(z)}{f(z)}$ has infinitely many
zeros. The proof of Theorem 4.2 is completed.
\vskip.2cm\par
{\bf Acknowledgements} The first author is supported by Guangdong National Natural Science Foundation(No.2014A030313422),and
the second author is supported by Training Plan Fund of Outstanding Young Teachers of
 Higher Learning Institutions of Guangdong Province of China(No.Yq20145084602) and Guangdong National Natural Science Foundation(No.2016A030313745)

   \vskip .5cm\par

\vskip.5cm\par\quad\\
Zhi-Bo Huang\\
School of Mathematical Sciences\\
South China Normal University\\
 Guangzhou, 510631, P.R.China\\
Email address:  huangzhibo@scnu.edu.cn
\vskip.5cm\par\quad\\
Ran-Ran Zhang\\
Department of Mathematics\\
Guangdong University of Education\\
Guangzhou, 510303, P.R.China\\
Email address:  zhangranran@gdei.edu.cn

\end{document}